\documentclass[a4,12pt]{article}

\usepackage{amssymb,amsmath}

\newtheorem{defn}{Definition}[section]

\newtheorem{exm}[defn]{Example}
\newtheorem{lemma}[defn]{Lemma}
\newtheorem{theorem}[defn]{Theorem}

\newtheorem{xdefn}{Definition.}
\newtheorem{xproposition}{Proposition.}
\newtheorem{xexm}{Example.}
\newtheorem{xlemma}{Lemma.}
\newtheorem{xtheorem}{Theorem.}
\newtheorem{xproof}{{\it Proof. }}

\newenvironment{definition}{\begin{defn}\em}{\end{defn}}
\newenvironment{example}{\begin{exm}\em}{\end{exm}}
\newenvironment{proof}{\begin{xproof}\em}{\end{xproof}}

\newenvironment{definition*}{\begin{xdefn}\em}{\end{xdefn}}
\newenvironment{example*}{\begin{xexm}\em}{\end{xexm}}
\newenvironment{proposition*}{\begin{xproposition}}{\end{xproposition}}
\newenvironment{lemma*}{\begin{xlemma}}{\end{xlemma}}
\newenvironment{theorem*}{\begin{xtheorem}}{\end{xtheorem}}

\def\qed{\hspace{0.3cm}{\rule{1ex}{2ex}}}
\newcommand\V{\bigvee}
\newcommand\Max{\mathrm{Max\,}}
\newcommand\ie{i.e.}

\newcommand\then{\&}
\newcommand\st{\mid}
\newcommand\ann{\mathrm{ann}}
\newcommand\Image{\mathrm{Im}}
\newcommand\Ann{\mathrm{Ann}}
\newcommand\End{\mathrm{Q}}
\newcommand\proj{\mathcal{P}}
\newcommand\B{\mathrm{B}}
\newcommand\hih{\mathcal{H}}
\newcommand\hik{\mathcal{K}}
\newcommand\caa{\mathfrak{A}}
\newcommand\cab{\mathfrak{B}}
\newcommand\ra{\rightarrow}
\newcommand\CC{\mathbb{C}}
\newcommand\upsegment{{\uparrow}}

\newcommand\ls{\mathrm{L}}

\begin{document}

\title{On quantales that classify C*-algebras\thanks{Research supported in part by FEDER and FCT/POCTI/POSI through the research units CAMGSD and CLC of IST and grant SFRH/BPD/11657/2002.}}
\author{David Kruml\thanks{On a post-doc leave from Masaryk University, Brno.} and Pedro Resende
\vspace*{2mm}\\ \small\it Departamento de Matem{\'a}tica, Instituto
Superior T{\'e}cnico, \vspace*{-2mm}\\ \small\it Av. Rovisco Pais,
1049-001 Lisboa, Portugal}

\date{~}

\maketitle

\vspace*{-1,3cm}
\begin{abstract}
The functor $\Max\!$  of Mulvey assigns to each unital C*-algebra $\caa$ the unital involutive quantale $\Max \caa$ of closed linear subspaces of $\caa$, and it has been remarked that it classifies unital C*-algebras up to $*$-isomorphism. In this paper we provide a proof of this and of the stronger fact that for every isomorphism $u:\Max\caa\ra\Max\cab$ of unital involutive quantales there is a $*$-isomorphism
$\widehat u:\caa\ra\cab$ such that $\Max\widehat u$ coincides with $u$ when restricted to the left-sided elements of $\Max\caa$. But we also show that isomorphisms $u:\Max\caa\ra\Max\cab$ may exist for which no isomorphism
$v:\caa\ra\cab$ is such that $\Max v=u$.
\vspace{0.1cm}
\vspace{0.1cm}
\end{abstract}

\section{Introduction}

This short paper is a followup to~\cite{KrPeReRo}, where various quantale~\cite{Mu86} based notions of spectrum of a C*-algebra were addressed from the point of view of their functorial properties. In particular, the functor $\Max\!$ from unital C*-algebras to unital involutive quantales, which was originally defined by Mulvey~\cite{Mu89} and was subsequently studied in~\cite{KrPeReRo,MuPe01,MuPe02} (see also the surveys~\cite{Mu02,PaRo00}) was seen to have no adjoints, therefore not providing the equivalence of categories that would be desired in order to consider $\Max\!$ a rightful ``spectrum functor''. However, as was also remarked in~\cite{KrPeReRo}, albeit without an explicit proof (one was presented in a talk~\cite{Re02:fields}), essentially from results of~\cite{GiKu71,MuPe01} it follows that $\Max\!$ classifies unital C*-algebras up to $*$-isomorphism, in the sense that any two unital C*-algebras $\caa$ and $\cab$ for which we have $\Max\caa\cong\Max\cab$ are necessarily $*$-isomorphic. Furthermore, from~\cite{Ro} it follows that $\Max\!$ is faithful~\cite{KrPeReRo}, thus leaving open the possibility that an equivalence of categories might be obtained between the category of unital C*-algebras and \emph{some} subcategory of the category of unital involutive quantales.

The aim of this paper is to provide some clarification regarding the above statements, and in it we prove the following result that, in particular, implies the classification theorem just mentioned:

\begin{theorem*}
Let $\caa$ and $\cab$ be unital C*-algebras, and let \[u:\Max\caa\ra\Max\cab\] be an isomorphism of unital involutive quantales. Then there is a $*$-isomorphism \[\widehat u:\caa\ra\cab\] such that $u(a)=\Max\widehat u(a)$ for every left-sided element $a\in\Max\caa$.
\end{theorem*}

In other words, $\Max\!$ is full on isomorphisms ``up to left-sided elements''.
However, we also show by means of an example that $\Max\!$ is not full on isomorphisms once the restriction to left-sided elements is dropped. The relevance of this observation follows from the following straightforward fact:

\begin{proposition*}
Let $C$ and $D$ be categories, and let $F:C\ra D$ be a functor. If $F$ is full on isomorphisms then its image $\Image (F)$ is a subcategory of $D$. If furthermore $F$ is faithful then $F:C\ra\Image( F)$ is an equivalence of categories.
\end{proposition*}

Hence, although $\Max\!$ has interesting properties, as discussed in~\cite{KrPeReRo}, it still does not provide us with an equivalence of categories in any obvious way.
In particular, our counterexample will show, for a particular C*-algebra $\caa$, that $\Max\caa$ has automorphisms which do not lie in the image of $\Max\!$, and thus, even though $\Max\!$ classifies unital C*-algebras up to $*$-isomorphism, it does not classify their automorphism groups.

For background on quantales and their modules we refer the reader to~\cite{Re:2forms}, whose notation and terminology we shall follow.

\section{Points of quantales}

Recall~\cite{Re:2forms} that an \emph{involutive} left module $M$ over a unital involutive quantale $Q$ is a left $Q$-module $M$ equipped with a \emph{symmetric sup-lattice 2-form} \[\varphi:M\times M\rightarrow 2\] (equivalently, an \emph{orthogonality relation} ${\perp}=\ker\varphi\subseteq M\times M$) satisfying, for all $a\in Q$ and $x,y\in M$,
\[\varphi(ax,y)=\varphi(x,a^*y)\;.\]
In addition, the \emph{annihilator} of an element $x\in M$ is the (left-sided) element
\[\ann(x)=\V\{a\in Q\st ax=0\}\;,\]
and $M$ is \emph{principal} if it has a \emph{generator}, by which we mean an element $x\in M$ such that $Qx=M$.

\begin{example}
Let $\caa$ be a unital C*-algebra, and let $\pi:\caa\ra\B(\hih)$ be a representation of $\caa$ on the Hilbert space $\hih$.
The sup-lattice $\proj(\hih)$ of norm-closed linear subspaces of $\hih$, with the usual orthogonality relation, is an involutive module over $\Max\caa$, with the action defined, for all $a\in\Max\caa$ and $P\in\proj(\hih)$, by
\[aP = \overline{\{\pi(A_1)(x_1)+\cdots+\pi(A_n)(x_n)\st A_i\in a,\ x_i\in P\}}\;,\]
where $\overline{(-)}$ denotes topological closure. This is the \emph{module induced by $\pi$}. Equivalently, we may view this module as a \emph{representation} as in~\cite{MuPe92,MuPe01,MuPe02,KrPeReRo}, \ie, the unital and involutive quantale homomorphism
\[\tilde\pi:\Max\caa\ra\End(\proj(\hih))\] given by $\tilde\pi(a)(P)=aP$, where $\End(\proj(\hih))$ is the quantale of endomorphisms of $\proj(\hih)$ with multiplication $f\then g=g\circ f$.
\end{example}

\begin{example}
Let $Q$ be a unital involutive quantale, and let $m\in Q$ be a left-sided element. Then the sup-lattice
\[\upsegment m=\{x\in Q\st m\le x\}\]
is an involutive left $Q$-module with the action and orthogonality relation being given by, for all $a\in Q$ and $x,y\in \upsegment m$,
\begin{eqnarray*}
ax &=& (a\then x)\vee m\;,\\
x\perp y &\iff& y^*\then x\le m\vee m^*\;.
\end{eqnarray*}
\end{example}

\begin{definition}
By a \emph{point} of a unital involutive quantale $Q$ will be meant (the isomorphism class of)
any principal involutive left $Q$-module $M$ for which there is a generator $x\in M$ such that $\ann(x)$ is a maximal left-sided element of $Q$.
\end{definition}

This notion of point differs from those of other papers~\cite{Kr:spatqu,KrPeReRo,MuPe01,MuRe} but it agrees with them insofar as irreducibility is concerned, since our points are necessarily irreducible representations~\cite[Th.\ 5.11]{Re:2forms}.

From~\cite{MuPe01} it follows that the unital involutive quantale $\Max\caa$ associated to a unital C*-algebra $\caa$ completely classifies the irreducible representations of $\caa$ up to unitary equivalence of representations. We shall now summarize these results. The first~\cite[Th.\ 9.1]{MuPe01} tells us that any point of $\Max\caa$ is induced by an irreducible representation of $\caa$. We state the aspects of this theorem which are important for us here, in the form presented in~\cite{Re:2forms}:

\begin{theorem}\label{thm:mp1}
    Let $\caa$ be a unital C*-algebra, and let $M$ be a point of
    $\Max \caa$.
    Then,
    \begin{enumerate}
    \item $M$ is induced by an 
    irreducible representation of $\caa$;
    \item $M$ is isomorphic as an involutive left $\Max\caa$-module to
    $\upsegment\ann(x)$
    for any generator $x$ of $M$.
    \end{enumerate}
\end{theorem}

From another result~\cite[Cor.\ 9.4]{MuPe01} it follows that also the relation of unitary equivalence of irreducible representations of a unital C*-algebra $\caa$ is determined by $\Max\caa$. We present here a different form of that result, along with a much shorter proof:

\begin{theorem}\label{thm:mp2}
    Let $\caa$ be a unital C*-algebra, and let
    $\pi_1:\caa\rightarrow\B(\hih_1)$ and $\pi_2:\caa\rightarrow\B(\hih_2)$
    be two irreducible representations of $\caa$ on Hilbert
    spaces $\hih_1$ and $\hih_2$. Then $\pi_1$ and $\pi_2$
    are unitarily equivalent if and only if the left $\Max \caa$-modules
    $\proj(\hih_1)$ and $\proj(\hih_2)$ which they induce are
    isomorphic.
\end{theorem}

\begin{proof}
    Let $f:\proj(\hih_1)\rightarrow \proj(\hih_2)$ be an isomorphism of left 
    $\Max \caa$-modules, and let $x\in \hih_1$. Writing
    $\bar x$ for the linear span of $x$, and $\Ann(x)\subseteq \caa$ for
    the annihilator of $x$ in $\caa$, we clearly have
    $\ann(\bar x)=\Ann(x)$. Also, $\ann(\bar x)=\ann(f(\bar x))$ because
    \[af(\bar x)=0\iff f(a\bar x)=0\iff a\bar x=0\]
    for all $a\in\Max \caa$. Finally, $f(\bar x)=\bar y$ for some $y\in \hih_2$,
    and $\Ann(y)=\ann(\bar y)$, and thus $\Ann(x)=\Ann(y)$. Hence,
    the two vectors $x$ and $y$ determine the same maximal
    left ideal
    of $\caa$, which means that the two representations $\pi_1$ and
    $\pi_2$ are determined by the same pure state of $\caa$,
    being thus equivalent.
    The converse, \ie, that equivalent representations determine
    isomorphic modules, is trivial. \qed
\end{proof}

We can indeed strengthen this result:

\begin{theorem}
Let $\caa$ be a unital C*-algebra, and let $M$ be a point
of $\Max\caa$, where $M$ equals $\proj(\hih)$ for some Hilbert space $\hih$ (and
the orthogonality relation is the usual one).
Then the left action of $\Max\caa$ on $M$ is induced by an irreducible representation
of $\caa$ on $\hih$.
\end{theorem}

\begin{proof}
From \ref{thm:mp1} it follows that there is an irreducible representation
of $\caa$ on a Hilbert space $\hik$
whose associated involutive left $\Max \caa$-module $\proj(\hik)$
is isomorphic to $\proj(\hih)$. It follows that $\hih$ and $\hik$ are isometrically
isomorphic because they have the same Hilbert dimension,
which coincides with the cardinality of any maximal pairwise orthogonal
set of atoms of $\proj(\hih)$, which of course is the same as the
cardinality of such a set taken from $\proj(\hik)$. Hence, the
irreducible representation on $\hik$ gives rise via the isomorphism to
an irreducible representation of $\caa$ on $\hih$, which furthermore
induces the original left $\Max \caa$-module structure of $\proj(\hih)$. \qed
\end{proof}

\section{Main results}

\begin{lemma}
Let $\caa$ and $\cab$ be unital C*-algebras, and let $u:\Max \caa\ra\Max \cab$ be an isomorphism of unital involutive quantales. Let also $\rho:\caa\ra\B(\hih)$ be an irreducible representation of $\caa$ on a Hilbert space $\hih$. Then there is an irreducible representation $\sigma:\cab\ra\B(\hih)$ of $\cab$ on $\hih$ such that $\tilde\rho=\tilde\sigma\circ u$ (where $\tilde\rho$ and $\tilde\sigma$ are the representations induced by
$\rho$ and $\sigma$, respectively).
\end{lemma}

\begin{proof}
It suffices to remark that $u$ obviously carries points to points because it is an isomorphism. In particular, $\tilde\rho\circ u^{-1}$ is a point of $\Max\cab$, and thus by the previous lemma there is an irreducible representation $\sigma:\cab\ra\B(\hih)$ such that $\tilde\sigma=\tilde\rho\circ u^{-1}$; equivalently, such that $\tilde\rho=\tilde\sigma\circ u$.
\qed
\end{proof}

\begin{lemma}
Let $\caa$ and $\cab$ be unital C*-algebras, and let $u:\Max \caa\ra\Max \cab$ be an isomorphism of unital involutive quantales. Let also $\rho_1:\caa\ra\B(\hih_1)$ 
and $\rho_2:\caa\ra\B(\hih_2)$ be irreducible representations of $\caa$, and
let $\sigma_1:\cab\ra\B(\hih_1)$ and $\sigma_2:\cab\ra\B(\hih_2)$ be irreducible
representations of $\cab$ such that $\tilde\sigma_i\circ u=\rho_i$ for $i=1,2$. Then
$\rho_1$ and $\rho_2$ are unitarily equivalent representations of $\caa$ if and only if $\sigma_1$ and $\sigma_2$ are unitarily equivalent representations of $\cab$.
\end{lemma}

\begin{proof}
First, $\rho_1$ and $\rho_2$ are equivalent if and only if $\tilde\rho_1$ and $\tilde\rho_2$ are equivalent representations of $\Max \caa$. Similarly, $\sigma_1$ and $\sigma_2$ are equivalent if and only if $\tilde\sigma_1$ and $\tilde\sigma_2$ are equivalent representations of $\Max\cab$. Finally, $u$ is an isomorphism and thus it preserves equivalence of representations, \ie, $\tilde\sigma_1$ and $\tilde\sigma_2$ are equivalent if and only if $\tilde\rho_1$ and $\tilde\rho_2$ are, since the latter equal $\tilde\sigma_1\circ u$
and $\tilde\sigma_2\circ u$, respectively. \qed
\end{proof}

\begin{theorem}
Let $\caa$ and $\cab$ be unital C*-algebras, and let \[u:\Max \caa\ra\Max \cab\] be an isomorphism of unital involutive quantales. Then there is a unital $*$-isomorphism $\widehat u:\caa\ra \cab$ such that $u$ coincides with $\Max{\widehat u}$ when restricted to left-sided elements.
\end{theorem}

\begin{proof}
Let $(\rho_i)_{i\in I}$ be a maximal family of pairwise inequivalent irreducible representations
of $\caa$ on Hilbert spaces $\hih_i$, with $i\in I$. By the
previous lemmas, there is a maximal family $(\sigma_i)_{i\in I}$ of pairwise inequivalent irreducible representations of $\cab$ on the same family of Hilbert spaces, such that for each $i\in I$ one has
\begin{equation}\label{eq:u}
\tilde\rho_i=\tilde\sigma_i\circ u\;.
\end{equation}
Hence, both $\caa$ and $\cab$
are isomorphic to weakly dense C*-subalgebras of the product of Von Neumann algebras
$\prod_{i\in I}\B(\hih_i)$ (which concretely consists of all the norm bounded $I$-indexed
families of operators).
More precisely, there is an embedding of unital C*-algebras
$\rho:\caa\ra\prod_{i\in I}\B(\hih_i)$ that to each $A\in\caa$ assigns the family
$(\rho_i(A))_{i\in I}\in\prod_{i\in I}\B(\hih_i)$, and another embedding
$\sigma:\cab\ra\prod_{i\in I}\B(\hih_i)$ that to each $A\in\caa$ assigns the family
$(\sigma_i(A))_{i\in I}$. Now let $P$ be a projection on $\hih_i$. Let us say that $P$ is an
\emph{open projection} (with respect to $\rho_i$) if $\ker P$ equals the annihilator in $\proj(\hih_i)$ of some $a\in\Max\caa$, in the following sense:
\[\ker P=\ann_{\rho_i}(a)=\V\{x\in\proj(\hih_i)\st \tilde\rho_i(a)(x)=0\}\;.\]
Similarly, let us call a projection $(P_i)_{i\in I}$ of $\prod_{i\in I}\B(\hih_i)$ \emph{open} with respect to
 $\rho$
if for each $i\in I$ the projection $P_i$ is open with respect to $\rho_i$. It turns out that a projection of
$\prod_{i\in I}\B(\hih_i)$ is open with respect to $\rho$ if and only if it is open with respect to $\sigma$,
because for each $i\in I$ we have
\[\ann_{\rho_i}(a)=\ann_{\sigma_i}(u(a))\;.\]
On the other hand, a projection is open with respect to $\rho$
if and only if it is an open $q$-set, in the sense of~\cite{GiKu71}, determined by the weakly dense inclusion of $\rho(\caa)$ into $\prod_{i\in I}\B(\hih_i)$ (that is, the support $e(K)$ of some subset $K\subseteq\rho(\caa)$).
Since the von Neumann algebra $\prod_{i\in I}\B(\hih_i)$ together with the open $q$-sets
determined in this way by the weakly dense inclusion of any unital C*-algebra in $\prod_{i\in I}\B(\hih_i)$ completely determine the C*-algebra as a C*-subalgebra of $\prod_{i\in I}\B(\hih_i)$~\cite[Th.\ 5.13]{GiKu71},
it follows that $\caa$ and $\cab$ are $*$-isomorphic. In particular, we have $\rho(\caa)=\sigma(\cab)$ and thus there is a $*$-isomorphism
\[\widehat u:\caa\ra\cab\]
defined by
\[\widehat u=(\sigma\vert_{\rho(\caa)})^{-1}\circ\rho\;.\]
Hence, we have $\rho=\sigma\circ \widehat u$, and thus also $\rho_i=\sigma_i\circ \widehat u$ for each $i\in I$.
From here it follows that
\begin{equation}\label{eq:whu}
\tilde\rho_i=\tilde\sigma_i\circ\Max{\widehat u}
\end{equation}
for each $i\in I$,
and thus we have two isomorphisms of unital involutive quantales,
\[u,\Max\widehat u:\Max\caa\ra\Max\cab\;,\]
satisfying similar conditions with respect to the points of $\Max\caa$ and $\Max\cab$, namely equations
(\ref{eq:u}) and (\ref{eq:whu}). This immediately implies that $u$ and $\Max\widehat u$ coincide on the left-sided elements of $\Max\caa$, because this quantale is known to be ``spatial on the left''~\cite{MuPe02} (equivalently, on the right), meaning precisely that its left-sided elements are separated by the points. \qed
\end{proof}

Hence, $\Max\!$ is full on isomorphisms ``up to left-sided elements''. However;

\begin{theorem}
$\Max\!$ is not full on isomorphisms.
\end{theorem}

\begin{proof}
Consider the commutative C*-algebra $\CC^2$. This has only two automorphisms, namely the identity and the map $(z,w)\mapsto (w,z)$, corresponding to the two permutations of the discrete two point spectrum of $\CC^2$.
Any automorphism of $\Max\CC^2$ is determined by its image on the atoms, which are the one dimensional subspaces of $\CC^2$. Consider then the following assignment to the atoms $\langle(z,w)\rangle$ of $\Max\CC^2$:
\begin{eqnarray*}
\langle(z,w)\rangle &\mapsto& \langle(w,z)\rangle\textrm{ if }z,w\neq 0\;,\\
\langle(z,0)\rangle &\mapsto& \langle(z,0)\rangle\;,\\
\langle(0,w)\rangle &\mapsto& \langle(0,w)\rangle\;.
\end{eqnarray*}
It is straightforward to check that this defines an automorphism of $\Max\CC^2$, which of course does not follow from any of the automorphisms of $\CC^2$. \qed
\end{proof}

In view of these results, one may be tempted to think that $\Max\caa$ has too much information in it and that attention should be focused on left-sided elements alone, since isomorphisms behave well with respect to these. A word of caution is in order, however, since previous studies of spectra based only on left-sided elements (equivalently, right-sided elements) have not been able to provide sufficiently powerful classification theorems: from~\cite{BoRoBo} it follows that the subquantale $\ls(\Max\caa)$ determines $\caa$ provided that we restrict to the class of post-liminary C*-algebras; and in~\cite{Ro} it is shown that the quantale $\ls(\Max\caa)$ equipped with the additional structure of a ``quantum frame'' is determined by the Jordan algebra structure of the self-adjoint elements of $\caa$.

Another natural way in which one may try to decrease the ``size'' of $\Max\caa$ is to take a quotient, instead of a subobject as just discussed. Observing that the good behavior of isomorphisms with respect to left-sided elements is a direct consequence of the spatiality of $\Max\caa$ ``on the left'', we may be led to replacing $\Max\caa$ by its ``spatial reflection'' in the hope that this will yield a better behaved functor. However, from~\cite{KrPeReRo} it follows that a functor does not arise in this way at all, because the spatialization of quantales with respect to their points is ill behaved, in particular not being a reflection.

\end{document}